\documentclass[10pt]{amsart}
\usepackage{hyperref}
\usepackage[bottom]{footmisc}
\usepackage{geometry}
\usepackage{xcolor}
\usepackage{graphicx}
\usepackage{enumitem}                               
\usepackage{amssymb}
\usepackage{mathrsfs}
\usepackage{amsthm}

\newtheorem*{theorem*}{Theorem}
\newtheorem*{lemma*}{Lemma}

\theoremstyle{definition}

\newtheoremstyle{definition*}
{\topsep}
{\topsep}
{}
{0pt}
{\bfseries}
{.}
{ }
{\thmname{#1}\thmnumber{ #2}\thmnote{ (#3)}}

\theoremstyle{definition*}
\newtheorem*{definition*}{Definition}
\newtheorem*{remark*}{Remark}
\newtheorem*{claim*}{Claim}
\newtheorem*{corollary*}{Corollary}

\begin{document}
	\baselineskip 13pt
	
	\begin{abstract} Suppose that $A$ is a finite nilpotent group of odd order acting good in the sense of \cite{EGJ} on the group $G$ of odd order. Under some additional assumptions we prove that the Fitting height of $G$ is bounded above by the sum of the numbers of primes dividing $|A|$ and $|C_G(A)|$ counted with multiplicities.
	\end{abstract}
	
	\title[good action of a nilpotent group]{good action of a nilpotent group\\ with regular orbits}
	\author{G\"{u}l\. In Ercan$^{*}$}
	\address{G\"{u}l\. In Ercan, Department of Mathematics, Middle East
		Technical University, Ankara, Turkey}
	\email{ercan@metu.edu.tr}
		\thanks{$^{*}$Corresponding author}
		\author{\.{I}sma\. Il \c{S}. G\"{u}lo\u{g}lu}
		\address{\.{I}sma\. Il \c{S}. G\"{u}lo\u{g}lu, Department of Mathematics, Do%
		\u{g}u\c{s} University, Istanbul, Turkey}
		\email{iguloglu@dogus.edu.tr}

	\maketitle
	
	\section{introduction} Let a finite group $A$ act on the group $G$ ``\textit{good}" as introduced in \cite{EGJ}, that is, $H=[H,B]C_H(B)$ for every subgroup $B$ of $A$ and for every $B$-invariant subgroup $H$ of $G$. It is expected that as a generalization of coprime action this concept may help to understand the real difficulties in studying a noncoprime action. In this framework \cite{EGJ} and \cite{EGJ2} are the first attempts generalizing some important coprime results to good action case. The present paper  partially extends one of the coprime  results of \cite{Tur4} to good action case from which the main result of \cite{EGJ} can be directly inferred.  Namely, we prove the following.\\
	
	\textbf{Theorem.} \textit{Let a finite nilpotent group $A$ of odd order which is $C_p\wr C_p$-free for any prime $p$ act good on the group $G$ of odd order. Assume further that $A$ acts with regular orbits, that is, for every subgroup $A_0$ of $A$ and every $A_0$-invariant irreducible
		elementary abelian section $S$ of G there is $v\in S$ with $C_A (v) = C_A (S)$. Let $B$ be a subgroup of $A$ such that $\bigcap_{a\in A}[G,B]^a=1.$ Then $h(G)\leq \ell(A:B)+ \ell(C_G(A)).$}\\
	
	\textbf{Corollary.} \textit{Let a finite nilpotent group $A$ of odd order which is $C_p\wr C_p$-free for any prime $p$ act good on the group $G$ of odd order. Assume further that $A$ acts with regular orbits, that is, for every subgroup $A_0$ of $A$ and every $A_0$-invariant irreducible
		elementary abelian section $S$ of G there is $v\in S$ with $C_A (v) = C_A (S)$. Then $h(G)\leq \ell(A)+ \ell(C_G(A)).$}\\
		
	Here, for a finite solvable group $S$, $h(S)$ is the Fitting height of $S$ and if $T\leq S$ $\ell(S:T)$ is the number of primes
	dividing the index $|S:T|$ counted with multiplicities.

	\section{proof of Theorem}
	Assume false. Choose a counterexample with minimum  $|A|+|A:B|$. If
	$A = B$ then $[G,A] = 1$ and the result holds, a contradiction. Hence $\ell(A:B)\geq 1$ and
	$h(G)\geq  2$. In view of \cite{EGJ}, $ G$ has an irreducible $A$-tower of
	height $h=h(G)$, that is a sequence of
	sections $P_1, \ldots, P_h$ of $G$ with $P_i = S_i / T_i$ where $S_i$ and $%
	T_i$ are $A$-invariant subgroups of $G$ satisfying conditions (1)-(8) of
	Definition 2.7 in \cite{EGJ}. It should be noted that we may assume that $T_h = 1$ and $
	S_h \leq F (G)$.
	
	We can assume that $P_1, \ldots, P_h$ is an irreducible $A$-tower of $G$ of height $h$
	with minimum $\prod_{i=1}^{h}|S_i|$. The minimality of $G$ yields $G = S_1 \ldots S_{h}$ and also $T_{h-1}=1.$
	To simplify the notation we set $V=S_h/\Phi(S_h)$ and $P=S_{h-1}$ and $X=PS_{h-2}S_{h-3}\ldots S_1.$ Note that $V$ is an irreducible $XA$-module. We shall proceed in a
	series of steps:\newline
	
	\textit{(1) $A$ acts faithfully on $G$, $A_1=C_A(P)\leq B$ and $(|P|,|A:B|)=1
		$. Furthermore we may assume that $\Phi(S_h)=1.$}
	
	\begin{proof}
		By induction applied to the action of $\bar{A}=A/Ker(A$ on $G)$ on $G$ with respect
		to the subgroup $\bar{B}$ we get $$h\leq \ell(A:BKer(A\; 
		\text{on}\;  G))+\ell(C_G(\bar{A}))\leq \ell(A:B)+\ell(C_G(A))$$ which is not the case. Therefore we may assume
		that $Ker(A$ on $G)=1$.
		
		We can observe that $A_1=C_A(P)$ centralizes all the subgroups $P,
		S_{h-2},\ldots ,S_1$ due to good action: Firstly we have $%
		[P_{h-2},A_1]=1$ by the three subgroups lemma, whence $
		[S_{h-2},A_1]=1$ by Proposition 2.2 (3) in \cite{EGJ}. Repeating the same argument we get
		the claim.
		
		Clearly $A_1\lhd A$. If $A_1\not \leq B$, by induction applied to the action
		of $A/A_1 $ on the group $PQ\ldots S_1$ with respect to the subgroup $%
		BA_1/A_1$ we have $$h-1\leq \ell(A/A_1:BA_1/A_1)+\ell(C_{PQ\ldots S_1}(A/A_1))\leq (\ell(A:B)-1)+\ell(C_{G}(A))$$ which is a contradiction.
		Thus $A_1\leq B$ and hence $(|P|,|A:B|)=1$ because $A_p$ centralizes $P$ by
		Proposition 2.5 in \cite{EGJ}. 
		
		Finally we may assume that $\Phi(S_h)=1$ as $h(G)=h(G/\Phi(S_h))$ holds.
	\end{proof}
	
	\textit{(2) For any subgroup $C$ of $A$ containing $B$ properly we have $P={%
			\langle [P,C]\rangle}^X.$ }
	
	\begin{proof}
		Set $P_0={\langle [P,C]\rangle}^X $, and $X_0=S_{h-2}\ldots S_1.$ Suppose
		that $P_0\ne P.$ Note that $P_0\lhd XC,$ and set $K=C_{X_0}(P/P_0).$ Then $%
		P_0K\lhd PX_0C\lhd XC$. Since $[P,C]\leq P_0$ we have $[X_0,C]\leq K$ by the
		three subgroups lemma. Then $[X,C]\leq P_0K$. Notice that $P_0S_{h-2}$ is
		normalized by $P_0K$. If $P\leq P_0K$ then $P$ normalizes $P_0S_{h-2}$ and
		so 
		\begin{equation*}
			P=[P,S_{h-2}]\leq [P,P_0S_{h-2}]\leq P_0S_{h-2}\cap P=P_0,
		\end{equation*}
		which is impossible. Thus we have $P\not \leq P_0K$ and so $P\not \leq
		\bigcap_{a\in A}[X,C]^a$. This forces that $P\cap \bigcap_{a\in
			A}[X,C]^a\leq \Phi(P)$ by condition (8) of Definition 2.7. Set $%
		Y=\bigcap_{a\in A}[X,C]^a$ and $\bar{X}=X/Y.$ An induction argument applied
		to the action of $A$ on $\bar{X}$ with respect to $C$ yields that 
		\begin{equation*}
			h(\bar{X})=h-1\leq \ell(A:C)+\ell(C_{\bar{X}}(A))=\ell(A:B)-\ell(C:B)+\ell(C_{\bar{X}}(A))
		\end{equation*}
		whence $h\leq \ell(A:B)+\ell(C_G(A)).$ This completes the proof of step \textit{(2)}.
	\end{proof}
	
	\textit{(3) $C_V(A)=0$.}
	
	\begin{proof}
		
		Suppose that $C_V(A)\ne 0$. Then $\ell(C_X(A)) \leq \ell(C_G(A))-1$. Since $\bigcap_{a\in A}[G,B]^a=1$, we have $h(G/[G,B])=h(G)=h$ and hence $h(G/[G,B]F(G))=h-1$. Thus $h(X/\bigcap_{a\in A}[X,B]^a)=h-1$. By induction applied to $X/\bigcap_{a\in A}[X,B]^a$ we get $$h-1\leq \ell(A:B)+ \ell(C_X(A))\leq \ell(A:B)+ \ell(C_G(A))-1$$ which is a contradiction.

	\end{proof}	
	\textit{(4) There exists $B_1\leq A$ such that $B\lhd B_1$; and an
		irreducible complex $XB_1$-submodule $M$ such that $M_{_{X}}$ is
		homogeneous, $P\not \leq Ker (X$ on $M), [X,B]\leq Ker (X$ on $M), C_M(B)=M$
		and $C_M(B_1)=0.$}
	
	\begin{proof}
		Clearly $V\not \leq [G,B]$ as $\bigcap_{a\in A}[G,B]^a=1.$ Note that $V_{XB}$
		is completely reducible as $XB\lhd \lhd XA$. Let $W$ be an irreducible $XB$%
		-submodule of $V$ such that $W\cap [G,B]=1$. Then $[W,B]=1$ and so $[X,B]\leq Ker(X$ on $W)$ by the three
		subgroups lemma. Therefore $W_{_{X}}$ is homogeneous. Let $U$ be the $X$%
		-homogeneous component of $V$ containing $W_{_{X}}$. Then $C_U(B)\ne 0,
		P\not \leq Ker(X$ on $U)$, and $[X,B]\leq Ker(X$ on $U)$. Set now $B_0=N_A(U)
		$. Then $U$ is an irreducible $XB_0$-module, and $B$ is properly contained
		in $B_0$ as $C_U(B_0)=0.$
		
		Let $\bar{k}$ be the algebraic closure of $k=\mathbb{F}_{p_{h}}$. Let $I$ be
		an irreducible submodule of $U\otimes _k \bar{k}.$ By the Fong-Swan theorem
		we may take an irreducible complex $XB_0$-module $M_0$ such that $Ker(X$ on $%
		M_0)=Ker(X$ on $U)$ and $M_0$ gives $I$ when reduced modulo $p_h$. Thus $%
		C_{M_0}(B)\ne 0,$ $P\not \leq Ker(X$ on $M_0)$, $[X,B]\leq Ker(X$ on $M_0)$
		and $C_{M_0}(B_0)=0.$ Observe that $B$ normalizes each $X$-homogeneous
		component of $M_0$ as $[X,B]\leq Ker(X$ on $M_0).$ Let now $M$ be an $X$%
		-homogeneous component of $M_0$ such that $C_{M}(B)\ne 0.$ Set $%
		B_1=N_{B_0}(M).$ Then we have $C_{M}(B_1)=0$ as $C_{M_0}(B_0)=0.$
		
		Suppose that $B$ is not normal in $B_1,$ and let $C=\langle B^{B_1} \rangle.$
		Then $[X,C]\leq Ker(X$ on $M)$. On the other hand, by \textit{(2)} we have $%
		P={\langle [P,C]\rangle}^X.$ This forces that $P\leq [X,C]\lhd Ker(X$ on $M)$%
		, which is not the case. Thus $B\lhd B_1$ whence $C_{M}(B)$ is a nonzero $%
		XB_1$-submodule of $M,$ and so $[M,B]=0$ by the irreducibility of $M$ as an $XB_1$-module. Then the claim holds.
	\end{proof}
	
	\textit{(4) Theorem follows.}
	
	\begin{proof}
		We consider the set of all pairs $(M_{\alpha}, C_{\alpha})$ such that $B\leq
		C_{\alpha}\leq B_1 $, $M_{\alpha}$ is an irreducible $XC_{\alpha}$-submodule
		of $M_{_{XC_{\alpha}}}$ and $C_{M_{\alpha}}(C_{\alpha})=0$. Choosing $(M_1,C)
		$ with $|C|$ minimum. Then $C_{M_1}(C_0)\ne 0$ for every $B\leq C_0<C$, $%
		(M_1)_{_{X}}$ is homogeneous and $Ker(X$ on $M_1)=Ker(X$ on $M).$
		
		Set now $\bar{X}=X/Ker(P$ on $M)$. We can observe that $[Z(\bar{P}),C]=1$:
		Otherwise, it follows by Theorem 3.3 in \cite{EGJ} that for any $\bar{P}$-homogeneous
		component $U$ of $(M_1)_{_{\bar{P}}}$, the module $U$ is $C$-invariant and $%
		\bar{X}=N_{\bar{X}}(U)C_{\bar{X}}(C)$. Then $C_{\bar{X}}(C)$ acts
		transitively on the set of all $\bar{P}$-homogeneous components of $M_1$.
		Clearly we have $[Z(\bar{P}),C]\leq Ker(\bar{P}$ on $U)$ and hence $[Z(\bar{P%
		}),C]=1$, as claimed. Thus if $\bar{P}$ is abelian, then $[P,C]\leq Ker(\bar{%
			P}$ on $M)$ and hence $P={\langle [P,C]\rangle}^X \leq Ker(\bar{P}$ on $M)$
		by \textit{(2)}, which is not the case. Therefore $\bar{P}$ is nonabelian.
		
		Let now $U$ be a homogeneous component of $(M_1)_{_{\Phi(\bar{P})}}$. Notice
		that $\Phi(\bar{P})\leq Z(\bar{P})$ by (5) of Definition 2.7 in \cite{EGJ} and so $[\Phi(%
		\bar{P}),C]=1$. Then $U$ is $C$-invariant. Set $\widehat{\bar{P}}=\bar{P}%
		/Ker({\bar{P}}$ on $U)$. Now $\Phi(\widehat{\bar{P}})=\widehat{\Phi(\bar{P})}
		$ is cyclic of prime order $p.$ Since $[Z(\bar{P}),C]=1$ we get $[X,C]\leq
		C_X(Z(\bar{P}))$ by the three subgroups lemma. Now clearly we have $%
		[X,C]\leq N_X(U).$ That is $X=N_X(U)C_X(C)$ as the action is good and so $%
		C_X(C)$ acts transitively on the set of all homogeneous components of $%
		(M_1)_{_{\Phi(\bar{P})}}$. Hence $M_1=\bigoplus_{t\in T}U^t$ where $T$ is a
		transversal for $N_X(U)$ in $X$ contained in $C_X(C).$ Notice that $N_{\bar{X%
			}C}(U)=N_{\bar{X}}(U)C.$ Set $X_1=C_X(\Phi(\bar{P}))$. Now $C_{XC}(\Phi(\bar{%
			P}))=X_1C\lhd XC$ and we have $[X,C]\leq X_1$. Then $X=X_1C_X(C).$ Clearly we have $PS_{h-2}\leq X_1\leq N_X(U)$ and $%
		X_1C\lhd XC\lhd \lhd XA.$ Recall that $P/\Phi(P)$ is an irreducible $XA$%
		-module and hence $P/\Phi(P)$ is completely reducible as an $X_1C$-module.
		Note that $\widehat{\bar{P}}/\Phi(\widehat{\bar{P}})\cong P/\Phi(P)C_P(U).$
		As $P/\Phi(P)$ is completely reducible we see that so is $P/\Phi(P)C_P(U)$.
		Hence $\widehat{\bar{P}}/\Phi(\widehat{\bar{P}})$ is also completely
		reducible.
		
		Since $\widehat{\Phi(\bar{P})}\leq \widehat{Z(\bar{P})}$, there
		is an $X_1C$-invariant subgroup $E$ containing $\widehat{\Phi(\bar{P})}$ so
		that 
		\begin{equation*}
			\widehat{\bar{P}}/\widehat{\Phi(\bar{P})}=\widehat{Z(\bar{P})}/\widehat{\Phi(%
				\bar{P})}\oplus E/\widehat{\Phi(\bar{P})}.
		\end{equation*}
		Then $\widehat{\bar{P}}=\widehat{Z(\bar{P})}E$ and hence $\widehat{Z(\bar{P})%
		}\cap E=Z(E).$ Clearly we have $({\widehat{\bar{P}}})^{\prime}=\widehat{\Phi(%
			\bar{P})}\leq Z(E).$ Also, 
		\begin{equation*}
			E/\widehat{\Phi(\bar{P})}\cap \widehat{Z(\bar{P})}/\widehat{\Phi(\bar{P})}=1
		\end{equation*}
		and hence $Z(E)\leq \widehat{\Phi(\bar{P})}$. Thus we have $Z(E)=\widehat{%
			\Phi(\bar{P})}=({\widehat{\bar{P}}})^{\prime}.$ As $E\unlhd \widehat{\bar{P}}$
		we get $\Phi(E)\leq \widehat{\Phi(\bar{P})}=Z(E).$ It follows that $%
		Z(E)=E^{\prime}=\Phi(E)=\widehat{\Phi(\bar{P})}$ is cyclic of prime order
		and hence $E$ is extraspecial. Now $[Z(\bar{P}),C]=1$ gives $[\widehat{Z(%
			\bar{P})},C]=1.$ Thus $[Z(E),C]=1.$
		
		Next we observe that $C_{C}(E)\leq B$: Otherwise there is a nonidentity
		element $b$ in $C\setminus B$ such that $[\widehat{\bar{P}},b]=1$ and hence $%
		[\bar{P},b]\leq Ker(\bar{P}$ on $U)$. Since $X=X_1C_X(C)\leq N_X(U)C_X(C)$
		we get $[\bar{P},b]\leq Ker(\bar{P}$ on $M)$. Set $C_1=B\langle b\rangle.$
		Recall that $[P,B]\leq Ker(P$ on $M)$ by \textit{(4)}. Then, by \textit{(2)}%
		, we have $P={\langle [P,C_1]\rangle}^X\leq Ker(P$ on $M)$ which is not the
		case. Therefore $C_{C}(E)\leq B$ as claimed.
		
		Notice that $p$ divides $|B/C_C(E)|$ if and only if $B_p\not \leq C_C(E)$
		which is impossible by Proposition 2.5 in \cite{EGJ} due to good action. This means by 
		\textit{(1)} that $p$ is coprime to $|C/C_C(E)|$. Note also that $B$ and
		hence $C_C(E)$ acts trivially on $U$ by \textit{(3)}. We apply now Lemma 2.1
		in \cite{Esp} to the action of the semidirect product $E(C/C_C(E))$ on the
		module $U$ and see that $C_U(C/C_C(E))\ne 0$. This final contradiction completes the proof.
		\end{proof}

\end{document}